\documentclass[fleqn,final]{zzz}
\usepackage{amssymb}
\usepackage{amsmath, amsfonts}
\usepackage{mathrsfs}
\usepackage{color}
\usepackage{tikz}
\usepackage{float}
\usepackage{graphicx}
\usepackage{rotating}
\usepackage[pdftex]{hyperref}
\usepackage{enumitem}
\usepackage{soul}
\usepackage{comment}
\usepackage{xparse}
\usepackage{chngcntr}
\usepackage{apptools}
\usepackage{caption}
\usepackage[normalem]{ulem}
\usepackage{subcaption}
\usepackage{arydshln}
\theoremstyle{lemma}

\theoremstyle{definition}
\newtheorem{ex}{Example}

\numberwithin{proposition}{section}
\numberwithin{remark}{section}
\numberwithin{equation}{section}
\numberwithin{theorem}{section}
\numberwithin{corollary}{section}

\DeclareMathOperator{\cosec}{csc}

\newcommand{\img}{\mathrm{i}}

\newcommand{\eul}{e}
\hypersetup{
colorlinks=true,
urlcolor=blue,
linkcolor=red,
citecolor=green
}

\NewDocumentCommand{\qfrac}{smm}{%
 \dfrac{\IfBooleanT{#1}{\vphantom{\big|}}#2}{\mathstrut #3}%
}

\makeatletter
\def\eqnarray{\stepcounter{equation}\let\@currentlabel=\theequation
\global\@eqnswtrue
\tabskip\@centering\let\\=\@eqncr
$$\halign to \displaywidth\bgroup\hfil\global\@eqcnt\z@
  $\displaystyle\tabskip\z@{##}$&\global\@eqcnt\@ne
  \hfil$\displaystyle{{}##{}}$\hfil
  &\global\@eqcnt\tw@ $\displaystyle{##}$\hfil
  \tabskip\@centering&\llap{##}\tabskip\z@\cr}

\def\endeqnarray{\@@eqncr\egroup
      \global\advance\c@equation\m@ne$$\global\@ignoretrue}

\def\@yeqncr{\@ifnextchar [{\@xeqncr}{\@xeqncr[5pt]}}
\makeatother

\newcommand{\C}{{\mathbb C}}

\newcommand{\Lp}{{\mathcal L}}

\newcommand{\Z}{{\mathbb Z}}

\definecolor{darkgreen}{rgb}{0.0, 0.21, 0.06}
\newcommand{\erf}{\mathrm{erf}}
\newcommand{\erfc}{\mathrm{erfc}}

\begin{document}

\renewcommand{\PaperNumber}{***}

\FirstPageHeading

\ShortArticleName{%
New identities for the Laplace transform and their applications
}

\ArticleName{%
New identities for the Laplace transform and their applications
}

\Author{Abdulhafeez A. Abdulsalam\,$^\dag\!\!\ $ and Ammar K. Mohammed\,$^\ast\!\!\ $}

\AuthorNameForHeading{A.~A.~Abdulsalam and A.~K.~Mohammed}
\Address{$^\dag$ Department of Mathematics, University of Ibadan, Ibadan, Oyo, Nigeria}
\EmailD{hafeez147258369@gmail.com, aabdulsalam030@stu.ui.edu.ng} 

\Address{$^\ast$ Department of Mathematics, University of Kirkuk, Kirkuk, Iraq}
\EmailD{scmma23043@uokirkuk.edu.iq} 


\ArticleDates{Received \today~in final form ????; 
Published online ????}

\Abstract{In this paper, we begin by applying the Laplace transform to derive closed forms for several challenging integrals that seem nearly impossible to evaluate. By utilizing the solution to the Pythagorean equation $a^2 + b^2 = c^2$, these closed forms become even more intriguing. This approach allows us to provide new integral representations for the error function. Additionally, by leveraging an identity we derived for the inverse Laplace transform and applying a result based on Srivastava and Y\"{u}rekli's identity, we provide a closed form for a nontrivial generalized integral.}



\Classification{11M06, 33B20, 42A38, 44A10, 44A35}
\allowdisplaybreaks
\section{Introduction}
The Gaussian integral $\int_{-\infty}^{\infty} \eul^{-x^2} \, \mathrm{d}x$ can be evaluated by a number of techniques. This results in the well-known identity $\Gamma\left(\frac{1}{2}\right) = \sqrt{\pi}$, and $\Gamma(z)$ is the gamma function \cite[(5.2.1)]{bib23}. The motivation of this paper arises from our observation that the Gaussian integral and several other integrals can be evaluated using the Laplace transform. A function $\phi(a, x)$ is said to be of exponential growth or order with respect to the variable $a$, if there exist constants $c$ and $M$ such that $|\phi(a, x)| \leq M \eul^{ca}$, where $\Re{a} > 0$. 
The Laplace transform of a piecewise function $\phi(a)$ of exponential order is defined by \cite[(1.14.17)]{bib23}
\[F(s) = \Lp\{\phi(a)\}(s) = \int_0^{\infty} \eul^{-as} \phi(a) \, \mathrm{d}a, \quad \Re{s} > 0,\]
where $\Re{z}$ and $\Im{z}$ denotes the real and imaginary part of the complex number $z$, respectively.
If $\phi(a)$ is continuous and $\phi'(a)$ is piecewise continuous on $[0, \infty)$, then \cite[(1.14.20)]{bib23}
\[\phi(a) = \Lp^{-1}\{F(s)\} = \frac{1}{2\pi \img} \lim_{T \to \infty} \int_{\sigma - \img T}^{\sigma + \img T} \eul^{as}  \Lp\{\phi(a)\}(s)\, \mathrm{d}s,
\]
where $\img  = \sqrt{-1}$ and $\sigma > c$. In this case, $\phi(a)$ is referred to as the inverse Laplace transform of $F(s)$. Notably, several integrals have been evaluated using the Laplace transform as demonstrated in \cite{gordon, apel}. In Section \ref{sec2.1}, we begin by utilizing the Laplace transform to derive closed forms for the following integrals, as well as several interesting integrals that can be reduced to or deduced from them:
\[\int_0^\infty \cos\left(bx \sqrt{x^2 + c^2}\right) \eul^{-ax^2} \, \mathrm{d}x, \quad \int_0^\infty \frac{\sin\left(bx \sqrt{x^2 + c^2}\right)}{x\sqrt{x^2 + c^2}} \eul^{-ax^2} \, \mathrm{d}x, \quad \int_0^\infty \frac{\sin\left(bx \sqrt{x + c}\right)}{\sqrt{x + c}} \eul^{-ax}\, \mathrm{d}x,
\]
for $\Re{a}, \Re{c} > 0$, $\Re{b} \geq 0$. If we let $\Phi(a, b, c, x)$, $\Psi(a, b, c, x)$, and $\varphi(a, b, c, x)$ represent the integrands of the three integrals respectively, then the relationships between the integrands are:
$$
\frac{\partial \Psi(a, b, c, x)}{\partial b} = \Phi(a, b, c, x), \quad x\Psi(a, bx, c, x) = \varphi(a, b, c^2, x^2).
$$
It is evident that these three integrals can be expressed in terms of the Laplace transform, with the first integral simplifying to the Gaussian integral when $b = 0$. Exploiting the fact that the $n$th partial derivatives of $\Phi(a, b, c, x)$, $\Psi(a, b, c, x)$, and $\varphi(a, b, c, x)$ are continuous with respect to $a$, $b$, $c$, and that their integrals are absolutely convergent over the domain of integration, we can deduce additional integrals once the initial three integrals have been evaluated. The list of integrals derivable from these three is extensive but not exhaustive. \\

\noindent
In general, consider a piecewise continuous function $\phi(a, x)$ on $[0, \infty)$ of exponential order with respect to the variable $a$, such that the limit $\lim_{x \to 0} \phi(a,x)$ exists. The methodology employed in Section \ref{sec2.1} can be applied to evaluate integrals of the form $\int_0^\infty \phi(a, x) \, \mathrm{d}x$. The error function and complementary error functions are defined, respectively, as \cite[(7.2.1), (7.2.2)]{bib23}
\[\erf{z} = \frac{2}{\sqrt{\pi}} \int_0^z \eul^{-u^2} \, \mathrm{d}u, \quad \erfc{z} = \frac{2}{\sqrt{\pi}} \int_z^\infty \eul^{-u^2} \, \mathrm{d}u = 1 - \erf{z}, \quad z \in \C.\]
From the first and second integrals, by applying the solution to the Pythagorean equation $a^2 + b^2 = c^2$, we derive several integral representations for the error function, one of which is:
\[
\erf{z}  = \frac{1}{\sqrt{2\pi}} \int_0^{2vz}  \frac{\sqrt{\sqrt{(v^2 - z^2)^2 + x^2} + (v^2 - z^2)}}{\sqrt{(v^2 - z^2)^2 + x^2}}\exp\left(\frac{1}{2} (v^2 - z^2) - \sqrt{(v^2 - z^2)^2 + x^2}\right) \, \mathrm{d}x,
\]
for $\Re{v} > 0$, $\Re(v^2) > \Re(z^2)$. The Fourier cosine and sine transforms are defined respectively by \cite[(1.14.9), (1.14.10)]{bib23}
\[\mathcal{F}_c\{f(a)\}(x) = \sqrt{\frac{2}{\pi}} \int_{0}^{\infty} f(a) \cos(ax) \, \mathrm{d}a, \quad \mathcal{F}_s\{f(a)\}(x) = \sqrt{\frac{2}{\pi}} \int_{0}^{\infty} f(a) \sin(ax) \, \mathrm{d}a.\]
The Mellin transform is defined by
\[\mathcal{M}\{f\}(s) = \int_0^\infty x^{s-1} f(x) \, \mathrm{d}x.\]
In Section \ref{sec2.2}, we derive more interesting identities using a formula due to Cauchy's residue theorem and an identity for the Fourier sine transform. Suppose $f(t)$ and $g(t)$ are absolutely and square integrable on $[0, \infty)$. Using the identity \cite[(1.14.14)]{bib23}
\[\int_0^{\infty} \mathcal{F}_s\{f(t)\}(x) \mathcal{F}_s\{g(t)\}(x)  \, \mathrm{d}x = \int_0^{\infty} f(t) g(t)\, \mathrm{d}t,
\]
and a result by Srivastava and Y\"{u}rekli \cite{stava}, we show that if $\Im(g(ix))$ is absolutely and square integrable on $[0, \infty)$, then 
\[\int_0^\infty \frac{\Im{g(\img x)}}{x^{\mu}} \, \mathrm{d}x = \frac{\cosec\left(\frac{\pi \mu}{2}\right)}{\Gamma(\mu)} \sqrt{\frac{\pi}{2}} \int_0^{\infty} x^{\mu -1} \mathcal{F}_s\{\Im{g(\img u)}\}(x) \, \mathrm{d}x, \quad \Re{\mu} > 0,\]
provided that each of the integrals exists. Applying this formula, we deduce the new generalized integral:
\[\int_0^\infty \frac{\sqrt{\sqrt{x^2+a^2} + a}}{\sqrt{x^2 + a^2}} x^{-\mu-1} \, \mathrm{d}x = -\frac{2^{\frac{1}{2}-2\mu} a^{-\frac{1}{2}-\mu} \Gamma(2\mu)}{\mu\Gamma(\mu)^2} \cosec\left(\frac{\pi \mu}{2}\right) \pi,\]
where $\Re{a} > 0$, $-\frac{1}{2} < \Re{\mu} < 0$. Throughout this paper, we verify our results using the Computer Algebra System (CAS) software \textsf{Mathematica 13}. All the formulas in Theorems \ref{thm1}--\ref{corfst} are new and have not been presented elsewhere in the literature.

\section{Results}
In this section, we present the main findings of this exposition, along with illustrative examples.
\subsection{Identities for the Laplace transform} \label{sec2.1}
In this subsection, we utilize the Laplace transform to prove several integral identities. Additionally, we provide examples to illustrate these identities. As an illustration, we begin by demonstrating the use of the Laplace transform in proving the closed form for the Gaussian integral $f(a) = \int_0^\infty \eul^{-ax^2} \, \mathrm{d}x$. By evaluating the Laplace transform of $f(a)$ and interchanging the order of integration, which is permissible by the absolute convergence of the integrals involved, we derive
\[F(s) = \int_0^\infty \int_0^\infty \eul^{-(s + x^2)a} \, \mathrm{d}a \, \mathrm{d}x = \int_0^\infty \frac{1}{x^2 + s} \, \mathrm{d}x = \frac{\pi}{2\sqrt{s}}.\]
The gamma function is defined by \cite[(5.2.1)]{bib23}
\[
\Gamma(p) = \int_0^\infty t^{p-1} \eul^{-t} \, \mathrm{d}t.
\]
Using the substitution $t = sa$, where $\Re(s) > 0$, we obtain
\[
\Gamma(p) = s^p \int_0^\infty a^{p-1} \eul^{-as} \, \mathrm{d}a = s^p \mathcal{L}\{a^{p-1}\}(s).
\]
This implies $\mathcal{L}^{-1}\left\{\frac{1}{s^p}\right\} = \frac{a^{p-1}}{\Gamma(p)}$. Replacing $s$ with $s + z$, where $\Re(s + z) > 0$, we obtain the more general form
\begin{equation} \label{hafpv1}
\mathcal{L}^{-1}\left\{\frac{1}{(s + z)^p}\right\}(a) = \frac{a^{p-1} \eul^{-az}}{\Gamma(p)}.
\end{equation}
Therefore, by evaluating the inverse Laplace transform of $F(s)$, we find $f(a) = \frac{\pi}{2}\mathcal{L}^{-1}\left\{\frac{1}{s^{\frac{1}{2}}}\right\}(a) = \frac{\pi a^{-\frac{1}{2}}}{2\Gamma\left(\frac{1}{2}\right)}$. Upon the substitution $x = t^{\frac{1}{2}}$, we find
\[f(a) = \frac{1}{2}\int_0^\infty t^{-\frac{1}{2}} \eul^{-at} \, \mathrm{d}t = \frac{\pi a^{-\frac{1}{2}}}{2\Gamma\left(\frac{1}{2}\right)}.\]
Setting $a=1$, we deduce $\Gamma\left(\frac{1}{2}\right) = \frac{\pi}{\Gamma\left(\frac{1}{2}\right)}$, and so $\Gamma\left(\frac{1}{2}\right)  = \sqrt{\pi}$. This implies $\int_{-\infty}^{\infty} \eul^{-ax^2} \, \mathrm{d}x = 2\int_0^{\infty} \eul^{-ax^2} \, \mathrm{d}x = \sqrt{\frac{\pi}{a}}$. It is well-known that Laplace evaluated the Gaussian integral in 1812. However, it is interesting to note that the above approach differs from the method employed by Laplace \cite[pp.~95, 96]{laplace} in evaluating this integral. It is unclear whether our application of the Laplace transform to evaluate the Gaussian integral is previously known. However, in the following theorems, we focus on providing closed forms for integrals that have not appeared in the literature before.
\begin{theorem}Let $a, b \in \mathbb{C}$, where $\Re{a} \geq 0$. Then \label{thm1}
\begin{equation}\label{eqthm1}
\int_0^\infty \cos\left(bx \sqrt{x^2 + 1}\right) \eul^{-ax^2} \, \mathrm{d}x =  \frac{\sqrt{\sqrt{a^2 + b^2} + a}}{2\sqrt{a^2 + b^2}} \eul^{\frac{a - \sqrt{a^2 + b^2}}{2}}\sqrt{\frac{\pi}{2}}.
\end{equation}
\end{theorem}

\begin{proof}
Consider the integral
\[f(a, b) =  \int_0^\infty \cos\left(abx \sqrt{x^2 + 1}\right) \eul^{-ax^2} \, \mathrm{d}x.\]
Evaluating the Laplace transform of $f(a, b)$ with respect to $a$, we have
\[
F(s) = \int_0^\infty \int_0^\infty \cos\left(abx \sqrt{x^2 + 1}\right) \eul^{-a(s + x^2)} \, \mathrm{d}x \, \mathrm{d}a.
\]
The intermediate steps for evaluating $F(s)$ can be found in Appendix \hyperref[appendix1]{\textcolor{red}{\textbf{A}}}. Finally, by applying \eqref{hafpv1} to determine the inverse Laplace transform of $F(s)$, we conclude
\begin{align*}
f(a, b) &= \frac{\pi\left(1 + b^2 + \sqrt{1 + b^2}\right)}{2\sqrt{2}(1+b^2)\sqrt{1 + \sqrt{1 + b^2}}} \mathcal{L}^{-1}\left\{\frac{1}{\left(s + \frac{b^2}{2\left(1 + \sqrt{1 + b^2}\right)}\right)^{\frac12}}\right\}(a)
\\&= \frac{\sqrt{\pi}\left(1 + b^2 + \sqrt{1 + b^2}\right)}{2\sqrt{2}\left(1+b^2\right) \sqrt{1 + \sqrt{1 + b^2}}} a^{-\frac{1}{2}} \eul^{-\frac{b^2 a}{2\left(1 + \sqrt{1+b^2}\right)}}.
\end{align*}
Replacing $b$ with $\frac{b}{a}$ in $f(a, b)$, we deduce
\begin{align*}
f\left(a, \frac{b}{a}\right) = \frac{1}{2}\sqrt{\frac{\pi}{2}} \frac{a^2 + b^2 + a \sqrt{a^2 + b^2}}{(a^2 + b^2)\sqrt{a + \sqrt{a^2 + b^2}}} = \frac{1}{2\sqrt{a^2 + b^2}} \frac{a + \sqrt{a^2 + b^2}}{\sqrt{a + \sqrt{a^2 + b^2}}} \sqrt{\frac{\pi}{2}}.
\end{align*}
Using $\frac{x}{\sqrt{x}} = \sqrt{x}$, we complete the proof of Theorem \ref{thm1}. 
\end{proof}
\begin{remark} The integral in Theorem \ref{thm1} can be succinctly expressed as $\frac{k}{l} e^{-m} \sqrt{\pi}$, where $k$, $l$ and $m$ are positive integers,  whenever $a$ and $b$ constitute the first two terms of a Pythagorean triple. More generally, for $a = p^2 - q^2$ and $b = 2pq$, where $p, q \in \mathbb{C}$ with $\Re(p^2) \geq \Re(q^2)$, the integral evaluates to $\frac{p}{2(p^2+q^2)} e^{-q^2} \sqrt{\pi}$.
\end{remark}
\begin{ex}The following result is valid:
\[\int_0^\infty \cos\left(4x \sqrt{x^2 + 1}\right) \eul^{-3x^2} \, \mathrm{d}x = \frac{\sqrt{\pi}}{5 \eul}.\]
\end{ex}
\noindent
Differentiating both sides of \eqref{eqthm1} with respect to $a$ once and twice yields the following corollaries, respectively.
\begin{corollary} For $\Re{a} > 0$, $b \in \mathbb{C}$. we have \label{corr1}
\begin{align*}
\int_0^\infty x^2 \cos\left(bx \sqrt{x^2 + 1}\right) \eul^{-ax^2} \, \mathrm{d}x &= \frac{\left(\left(a^2 - b^2\right) + (a - b^2) \sqrt{a^2 + b^2}\right)\sqrt{\sqrt{a^2 + b^2} - a}}{4b(a^2 + b^2)^{\frac{3}{2}}} 
\\&\qquad \qquad \times \eul^{\frac{a - \sqrt{a^2 + b^2}}{2}} \sqrt{\frac{\pi}{2}}.
\end{align*} 
\end{corollary}
\begin{ex}The following result is valid:
\[
\int_0^\infty x^2 \cos\left(4x \sqrt{x^2 + 1}\right) \eul^{-3x^2} \, \mathrm{d}x = -\frac{9 \sqrt{\pi }}{250 \eul}.\]
\end{ex}
\begin{corollary} For $\Re{a} > 0$, $b \in \mathbb{C}$. we have \label{corr2}
\begin{align*}
\int_0^\infty &x^4  \cos\left(bx \sqrt{x^2 + 1}\right) \eul^{-ax^2} \, \mathrm{d}x 
\\&\quad= \frac{\left(\sqrt{a^2+b^2}-a\right)^{\frac{3}{2}}}{8 b^3 \left(a^2+b^2\right)^{\frac{5}{2}}}\left(6 a^4-9 a^2 b^2-6 a^3 b^2-3 b^4-6 a b^4 +a^2 b^4+b^6  \right.
\\& \qquad \qquad \qquad \qquad \qquad \qquad   \left. +  \left(6a^3 - 12ab^2 - 6a^2b^2\right)\sqrt{a^2+b^2}\right) \eul^{\frac{a - \sqrt{a^2 + b^2}}{2}} \sqrt{\frac{\pi}{2}}.
\end{align*}
\begin{ex}The following result is valid:
\[\int_0^\infty x^4 \cos\left(4x \sqrt{x^2 + 1}\right) \eul^{-3x^2} \, \mathrm{d}x = -\frac{137 \sqrt{\pi }}{12500 \eul}.\]
\end{ex}
\end{corollary}
\noindent
By making necessary substitutions and differentiating both sides of  \eqref{eqthm1} $n$ times with repsect to $a$, we obtain the following result.
\begin{corollary}For $\Re{a}, \Re{c} > 0$, $b \in \mathbb{C}$, and $n$ a nonnegative integer, we have \label{cor2}
\begin{align*}
&\int_0^\infty x^{2n} \cos\left(bx \sqrt{x^2 + c^2}\right) \eul^{-ax^2} \, \mathrm{d}x 
\\&\qquad \qquad \qquad= (-1)^n \frac{1}{2}\sqrt{\frac{\pi}{2}} \frac{\partial^n}{\partial a^n} \left(\frac{c \sqrt{\sqrt{a^2 c^4 + b^2 c^4} + a c^2}}{\sqrt{a^2 c^4 + b^2 c^4}} \eul^{\frac{a c^2 - \sqrt{a^2 c^4 + b^2 c^4}}{2}}\right).
\end{align*}
\end{corollary}
\noindent
Differentiating both sides of \eqref{eqthm1} with respect to $b$ once and twice, we obtain the following corollaries, respectively.
\begin{corollary} For $\Re{a} > 0$, $b \in \mathbb{C}$. we have \label{corr3}
\begin{align*}
&\int_0^\infty x \sqrt{x^2 + 1} \sin\left(bx \sqrt{x^2 + 1}\right) \eul^{-ax^2} \, \mathrm{d}x 
\\&\qquad \qquad \qquad= \frac{\left(a^2 + 2a + b^2 + (a + 1)\sqrt{a^2 + b^2}\right)\sqrt{\sqrt{a^2 + b^2} - a}}{4(a^2 + b^2)^{\frac{3}{2}}}\eul^{\frac{a - \sqrt{a^2 + b^2}}{2}} \sqrt{\frac{\pi}{2}}.
\end{align*} 
\end{corollary}
\begin{ex}The following result is valid:
\[\int_0^\infty x \sqrt{x^2 + 1} \sin\left(4x \sqrt{x^2 + 1}\right) \eul^{-3x^2} \, \mathrm{d}x =  \frac{51 \sqrt{\pi }}{500 \eul}.\]
\end{ex}
\begin{corollary} For $\Re{a} > 0$, $b \in \mathbb{C}$. we have \label{corr4}
\begin{align*}
\int_0^\infty& x^2 (x^2 + 1) \cos\left(bx \sqrt{x^2 + 1}\right) \eul^{-ax^2} \, \mathrm{d}x 
\\&\qquad= \frac{\sqrt{\sqrt{a^2+b^2}+a}}{8 \left(a^2+b^2\right)^{\frac{5}{2}}}\left(9 a^2-2 a^3-3 b^2-2 a b^2 -a^2 b^2 - b^4 \right.
\\& \qquad \qquad \qquad \qquad \qquad \quad  \left.+ (4a^2 - 6a - 2b^2)\sqrt{a^2 + b^2}\right)  \eul^{\frac{a - \sqrt{a^2 + b^2}}{2}} \sqrt{\frac{\pi}{2}}.
\end{align*}
\end{corollary}
\begin{ex}The following result is valid:
\[\int_0^\infty x^2 (x^2 + 1) \cos\left(4x \sqrt{x^2 + 1}\right) \eul^{-3x^2} \, \mathrm{d}x =  -\frac{587 \sqrt{\pi }}{12500 \eul}.\]
\end{ex}
\noindent
By making necessary substitutions and differentiating both sides of \eqref{eqthm1} $n$ times with respect to $b$, we obtain the following result.
\begin{corollary}For $\Re{a}, \Re{c} > 0$, $b \in \mathbb{C}$, and $n$ a nonnegative integer, we have  \label{cor3}
\begin{align*}
&\int_0^\infty x^n (x^2 + c^2)^{\frac{n}{2}} \cos\left(bx \sqrt{x^2 + c^2} + \frac{n\pi}{2}\right) \eul^{-ax^2} \, \mathrm{d}x 
\\&\qquad \qquad \qquad = \frac{1}{2}\sqrt{\frac{\pi}{2}} \frac{\partial^n}{\partial b^n} \left(\frac{c \sqrt{\sqrt{a^2 c^4 + b^2 c^4} + a c^2}}{\sqrt{a^2 c^4 + b^2 c^4}} \eul^{\frac{a c^2 - \sqrt{a^2 c^4 + b^2 c^4}}{2}}\right).
\end{align*}
\end{corollary}

\begin{theorem}Let $a, b, c \in \C$, where $\Re{a}, \Re{b} > 0$. Then \label{thm2}
\begin{equation}\label{yind1}
\int_0^\infty \frac{\sqrt{\sqrt{x^2 + a^2} + a}}{\sqrt{x^2 + a^2}} \cos(c x) \eul^{-b\sqrt{x^2 + a^2}} \, \mathrm{d}x = \frac{\sqrt{b + \sqrt{b^2 + c^2}}}{\sqrt{b^2 + c^2}} \eul^{-a\sqrt{b^2 + c^2}} \sqrt{\frac{\pi}{2}}.
\end{equation}
\end{theorem}

\begin{proof}
Let us denote
\[F(a, b, c) = \int_0^\infty \frac{\sqrt{\sqrt{x^2 + a^2} + a}}{\sqrt{x^2 + a^2}} \cos(c x) \eul^{-b\sqrt{x^2 + a^2}} \, \mathrm{d}x.\]
By substituting $x = au$, we obtain
\[F(a, b, c) = \sqrt{a} \int_0^\infty \frac{\sqrt{\sqrt{u^2 + 1} + 1}}{\sqrt{u^2 + 1}} \cos(acu) \eul^{-ab\sqrt{u^2 + 1}} \, \mathrm{d}u.\]
Further substitution of $u = \frac{\sqrt{1 -v^2}}{v}$ gives
\[F(a, b, c) = \sqrt{a} \int_0^1 v^{-\frac{3}{2}} (1- v)^{-\frac{1}{2}} \eul^{-\frac{ab}{v}} \cos\left(\frac{a c \sqrt{1-v^2}}{v}\right) \mathrm{d}v.\]
Finally, substituting $v = \frac{1}{2w^2 + 1}$ yields
\begin{equation}\label{leqn2}
F(a, b, c) = 2\sqrt{2a}\eul^{-ab}\int_0^\infty \cos\left(2acw\sqrt{w^2 +1}\right)\eul^{-2abw^2}\, \mathrm{d}w.
\end{equation}
Utilizing Theorem \ref{thm1} in \eqref{leqn2}, we conclude the proof of Theorem \ref{thm2}.
\end{proof}

\begin{ex}The following results are valid: \label{examp6}
\begin{align*}
\int_0^\infty \frac{\sqrt{\sqrt{x^2 + 4} + 2}}{\sqrt{x^2 + 4}} \cos(4x) \eul^{-3\sqrt{x^2 + 4}} \, \mathrm{d}x  &= \frac{2 \sqrt{\pi }}{5 e^{10}}, \\
\int_0^\infty \frac{\sqrt{\sqrt{x^2 + 9} + 3}}{\sqrt{x^2 + 9}} \cos(12 x) \eul^{-5\sqrt{x^2 + 9}} \, \mathrm{d}x &= \frac{3 \sqrt{\pi }}{13 e^{39}}.
\end{align*}
\end{ex}

\noindent
We note that the first two closed forms in Example \ref{examp6} admit the form $\frac{a}{b} \eul^{-ab}$. This observation motivates the following corollary.
\begin{corollary} Let $\alpha, \beta \in \C$, where $\Re(\alpha) > 0$ and $\Re(\alpha^2) > \Re(\beta^2)$. Then \label{corf}
\begin{equation} \label{ne1}
\int_0^\infty \frac{\sqrt{\sqrt{x^2 + \alpha^2} + \alpha}}{\sqrt{x^2 + \alpha^2}} \cos(2\alpha \beta x) \eul^{-(\alpha^2 - \beta^2)\sqrt{x^2 + \alpha^2}} \, \mathrm{d}x = \frac{\alpha}{\alpha^2 + \beta^2} \eul^{-\alpha(\alpha^2 + \beta^2)} \sqrt{\pi}.
\end{equation}
\end{corollary}

\begin{proof}
Utilizing the fact that integer solutions for the Pythagorean equation $b^2 + c^2 = d^2$ are of the form $b = \alpha^2 - \beta^2$, $c = 2\alpha \beta$, and $d = \alpha^2 + \beta^2$, where $\alpha, \beta \in \Z^{+}$, $\alpha > \beta$, we substitute these values alongside $a=\alpha$ into Theorem \ref{thm2} to conclude the proof of \eqref{ne1}. Regarding the conditions, for the integral to converge, we require $\Re(\alpha^2 - \beta^2) =\Re(\alpha^2) - \Re(\beta^2) > 0$.  This completes the proof of Corollary \ref{corf}.
\end{proof}

\noindent
By differentiating both sides of  \eqref{yind1} with repsect to $b$, we obtain the following result.
\begin{corollary} Let $a, b, c \in \C$, where $\Re{a}, \Re{b} > 0$. Then \label{corfrmthm2}
\begin{align*}
\int_0^\infty &\sqrt{\sqrt{x^2 + a^2} + a} \cos(c x) \eul^{-b\sqrt{x^2 + a^2}} \, \mathrm{d}x 
\\&\qquad= \frac{\left((2ab -1)\sqrt{b^2+c^2} + 2b\right)\sqrt{\sqrt{b^2+c^2}+b}}{2 \left(b^2+c^2\right)^{\frac{3}{2}}} \eul^{-a \sqrt{b^2+c^2}} \sqrt{\frac{\pi }{2}} .
\end{align*}
\end{corollary}

\begin{ex}The following result is valid: 
\[\int_0^\infty \sqrt{\sqrt{x^2 + 1} + 1}\, \cos(4x) \eul^{-3\sqrt{x^2 + 1}} \, \mathrm{d}x  = \frac{31\sqrt{\pi}}{125\eul^{5}}.\]
\end{ex}

\noindent
By differentiating both sides of  \eqref{yind1} twice with repsect to $b$, we obtain the following result.
\begin{corollary} Let $a, b, c \in \C$, where $\Re{a}, \Re{b} > 0$. Then \label{cor2fromthm2}
\begin{align*}
\int_0^\infty &\sqrt{\sqrt{x^2 + a^2} + a} \sqrt{x^2 + a^2}  \cos(c x) \eul^{-b\sqrt{x^2 + a^2}} \, \mathrm{d}x 
\\&\qquad= \frac{\sqrt{\sqrt{b^2+c^2}+b}}{4 \left(b^2+c^2\right)^{\frac{5}{2}}} \left(9 b^2 - 4 a b^3 + 4 a^2 b^4 - 3 c^2 - 4 a b c^2 + 4 a^2 b^2 c^2\right.
\\&\qquad \qquad \qquad \qquad \qquad \quad \left.  + (8 a b^2 - 6 b - 4 a c^2 ) \sqrt{b^2 + c^2}\right) \eul^{-a \sqrt{b^2+c^2}} \sqrt{\frac{\pi }{2}}.
\end{align*}
\end{corollary}

\begin{ex}The following result is valid: 
\[\int_0^\infty \sqrt{\sqrt{x^2 + 1} + 1}\sqrt{x^2 + 1}\, \cos(4x) \eul^{-3\sqrt{x^2 + 1}} \, \mathrm{d}x  = \frac{583\sqrt{\pi}}{6250\eul^{5}}.\]
\end{ex}

\noindent
By differentiating both sides of  \eqref{yind1} with repsect to $c$, we obtain the following result.
\begin{corollary}Let $a, b, c \in \C$, where $\Re{a}, \Re{b} > 0$, $\Re{c} \geq 0$. Then \label{yinde1}
\begin{align*}
\int_0^\infty &\frac{x\sqrt{\sqrt{x^2 + a^2} + a}}{\sqrt{x^2 + a^2}}  \sin(c x) \eul^{-b\sqrt{x^2 + a^2}} \, \mathrm{d}x 
\\&\quad= \frac{\left((2ab + 1)\sqrt{b^2 + c^2} + 2 \left(b + a b^2 + a c^2\right)\right)\sqrt{\sqrt{b^2+c^2}-b}}{2\left(b^2+c^2\right)^{\frac{3}{2}}} \eul^{-a \sqrt{b^2+c^2}} \sqrt{\frac{\pi}{2}}.
\end{align*}
\end{corollary}

\begin{ex}The following result is valid: 
\[
\int_0^\infty \frac{x\sqrt{\sqrt{x^2 + 1} + 1}}{\sqrt{x^2 + 1}} \sin(4x) \eul^{-3\sqrt{x^2 + 1}} \, \mathrm{d}x  = \frac{91\sqrt{\pi}}{250\eul^{5}}.\]
\end{ex}

\noindent
By differentiating both sides of  \eqref{yind1} twice with repsect to $c$, we obtain the following result.
\begin{corollary}Let $a, b, c \in \C$, where $\Re{a}, \Re{b} > 0$. Then \label{yinde2}
\begin{align*}
\int_0^\infty &\frac{x^2\sqrt{\sqrt{x^2 + a^2} - a}}{\sqrt{x^2 + a^2}}  \cos(c x) \eul^{-b\sqrt{x^2 + a^2}} \, \mathrm{d}x 
\\&\qquad= \frac{\sqrt{\sqrt{b^2+c^2}+b}}{4 \left(b^2+c^2\right)^{\frac{5}{2}}} \left(9 b^2 - 4 a b^3 - 3 c^2 - 4 a b c^2 - 4 a^2 b^2 c^2 - 4 a^2 c^4\right.
\\&\qquad \qquad \qquad \qquad \qquad \quad \left.  + (8 a b^2 - 6 b - 4 a c^2 ) \sqrt{b^2 + c^2}\right) \eul^{-a \sqrt{b^2+c^2}} \sqrt{\frac{\pi }{2}}.
\end{align*}
\end{corollary}

\begin{ex}The following result is valid: 
\[
\int_0^\infty \frac{x^2\sqrt{\sqrt{x^2 + 1} + 1}}{\sqrt{x^2 + 1}} \cos(4x) \eul^{-3\sqrt{x^2 + 1}} \, \mathrm{d}x  = -\frac{1917\sqrt{\pi}}{6250\eul^{5}}.\]
\end{ex}

\noindent
By differentiating both sides of  \eqref{yind1} with repsect to $a$, we obtain the following result.
\begin{corollary}Let $a, b, c \in \C$, where $\Re{a}, \Re{b} > 0$. Then \label{yindef5}
\begin{align*}
\int_0^\infty &\frac{\sqrt{\sqrt{x^2+a^2}+a} \left((2ab-1) \sqrt{x^2+a^2} + 2a\right) \cos(cx)}{\left(x^2+a^2\right)^{\frac{3}{2}}} \eul^{-b \sqrt{x^2+a^2}} \, \mathrm{d}x 
\\&\qquad\quad= \sqrt{\sqrt{b^2 + c^2} +  b} \eul^{-a \sqrt{b^2+c^2}} \sqrt{2\pi}.
\end{align*}
\end{corollary}

\begin{ex}The following result is valid: 
\[
\int_0^\infty \frac{\sqrt{\sqrt{x^2+1}+1} \left(5\sqrt{x^2+ 1} + 2\right) \cos(4x)}{\left(x^2+1\right)^{\frac{3}{2}}} \eul^{-3\sqrt{x^2+1}} \, \mathrm{d}x  = \frac{4\sqrt{\pi}}{\eul^{5}}.\]
\end{ex}

\noindent
By differentiating both sides of  \eqref{yind1} twice with repsect to $a$, we obtain the following result.
\begin{corollary}Let $a, b, c \in \C$, where $\Re{a}, \Re{b} > 0$. Then \label{yindef6}
\begin{align*}
\int_0^\infty &\frac{\sqrt{\sqrt{x^2+a^2}+a}}{\left(x^2+a^2\right)^{\frac{5}{2}}} \left(9 a^2 - 4 a^3 b + 4 a^4 b^2 - 3 x^2 - 4 a b x^2 + 4 a^2 b^2 x^2 \right.
\\& \qquad \qquad \qquad \qquad \left.+ \left(8a^2b - 6a - 4bx^2\right) \sqrt{x^2+a^2}\right) \cos(cx) \eul^{-b \sqrt{x^2+a^2}} \, \mathrm{d}x 
\\&\qquad= 2\sqrt{\sqrt{b^2 + c^2} + b} \sqrt{b^2+c^2} \eul^{-a \sqrt{b^2+c^2}} \sqrt{2\pi}.
\end{align*}
\end{corollary}

\begin{ex}The following result is valid: 
\[\int_0^\infty \frac{\sqrt{\sqrt{x^2+1}+1} \left(\left(18-12 x^2\right) \sqrt{x^2+1} + 33 + 21 x^2\right) \cos(4x)}{\left(x^2+1\right)^{\frac{5}{2}}} \eul^{-3\sqrt{x^2+1}} \, \mathrm{d}x  = \frac{40\sqrt{\pi}}{\eul^{5}}.\]
\end{ex}

\begin{theorem}Let $a, b \in \C$, where $\Re{a},\, \Re{b} \geq 0$. Then \label{thm3}
\begin{equation}\label{eqthm2}
\int_0^\infty \frac{\sin\left(bx \sqrt{x^2 + 1}\right)}{x\sqrt{x^2 + 1}} \eul^{-ax^2} \, \mathrm{d}x = \frac{\pi}{2}\erf\left(\sqrt{\frac{\sqrt{a^2 + b^2} - a}{2}}\right).
\end{equation}
\end{theorem}

\begin{proof}
Consider the integral
\[f(a, b) = \int_0^\infty \frac{\sin\left(abx \sqrt{x^2 + 1}\right)}{x\sqrt{x^2 + 1}} \eul^{-ax^2} \, \mathrm{d}x.\]
Evaluating the Laplace transform of $f(a, b)$ with respect to $a$, we have
\[
F(s) = \int_0^\infty \int_0^\infty \frac{\sin\left(abx \sqrt{x^2 + 1}\right)}{x\sqrt{x^2 + 1}} \, \mathrm{d}x \, \mathrm{d}a.
\]
The intermediate steps for evaluating $F(s)$ can be found in Appendix \hyperref[appendix2]{\textcolor{red}{\textbf{B}}}. Upon rationalization, we arrive at
\[
F(s) = \frac{1}{\sqrt{b^2 - 4s^2 + 4s}}\left(\frac{\pi \sqrt{b^2+2s - b\sqrt{b^2-4 s^2+4 s}}}{2\sqrt{2}s} - \frac{\pi \sqrt{b^2+2s + b\sqrt{b^2-4 s^2+4 s}}}{2\sqrt{2}s}\right).
\]
Using the identity $\sqrt{x} -\sqrt{y} = \sqrt{x + y - 2\sqrt{xy}}$, we find
\[
F(s) = \frac{\pi}{2\sqrt{2}s\sqrt{b^2 - 4s^2 + 4s}} \cdot \sqrt{2} \sqrt{b^2 + 2s - 2s \sqrt{1 + b^2}} = \frac{\sqrt{b^2 + 2s - 2s \sqrt{1 + b^2}}}{2s \sqrt{b^2 - 4s^2 + 4s}} \pi.
\]
Rationalizing and simplifying, we obtain
\[
F(s) = \frac{\sqrt{(b^2 + 2s)^2 - 4s^2(1 + b^2)}}{2s \sqrt{b^2 - 4s^2 + 4s} \sqrt{b^2 + 2s + 2s\sqrt{1+b^2}}} \pi  =\frac{b}{2s\sqrt{b^2 + 2s + 2s\sqrt{1+b^2}}}\pi.
\]
According to the convolution theorem \cite[(1.14.31)]{bib23},
\[ 
\mathcal{L}^{-1}\{F(s) \cdot G(s)\}(a) = \mathcal{L}^{-1}\{F(s)\}(a) \ast \mathcal{L}^{-1}\{G(s)\}(a),
\]
where the convolution is defined by
\[\mathcal{L}^{-1}\{F(s)\}(a) \ast \mathcal{L}^{-1}\{G(s)\}(a) = \int_0^a \mathcal{L}^{-1}\{F(s)\}(a-u) \mathcal{L}^{-1}\{G(s)\}(u) \, \mathrm{d}u.
\]
Evaluating the inverse Laplace transform of $F(s)$, we have
\[f(a, b) = \pi \mathcal{L}^{-1}\left\{\frac{1}{s} \cdot \frac{b}{2\sqrt{b^2 + 2s + 2s\sqrt{1+b^2}}}\right\}(a).
\]
Applying the convolution theorem, we derive
\begin{align*}
f(a, b) &= \pi \left(\mathcal{L}^{-1}\left\{\frac{1}{s}\right\}(a) \ast \mathcal{L}^{-1}\left\{\frac{b}{2\sqrt{b^2 + 2s + 2s\sqrt{1+b^2}}}\right\}(a) \right)
\\&=  \pi \left(1 \ast \mathcal{L}^{-1}\left\{\frac{b}{2\sqrt{b^2 + 2s + 2s\sqrt{1+b^2}}}\right\}(a)\right)
\\&=  \frac{b \pi}{2\sqrt{2}\left(1 + \sqrt{1+b^2}\right)} \int_0^a \mathcal{L}^{-1}\left\{ \frac{1}{\left(s + \frac{b^2}{2\left(1 + \sqrt{1+b^2}\right)}\right)^{\frac{1}{2}}}\right\}(u) \, \mathrm{d}u.
\end{align*}
Finally, using \eqref{hafpv1}, we conclude that
\[
f(a, b) = \frac{b}{2\sqrt{1 + \sqrt{b^2 + 1}}} \sqrt{\frac{\pi}{2}} \int_0^a u^{-\frac{1}{2}} \eul^{-\frac{b^2 u}{2 \left(\sqrt{b^2+1}+1\right)}} \, \mathrm{d}u = \frac{\pi}{2}\erf\left(\frac{b\sqrt{a}}{\sqrt{2\left(1 + \sqrt{1+b^2}\right)}}\right).\]
Replacing $b$ with $\frac{b}{a}$ in $f(a, b)$, we deduce
\[
f\left(a, \frac{b}{a}\right) = \frac{\pi}{2}\erf\left(\frac{b}{\sqrt{2\left(\sqrt{a^2 +b^2} + a\right)}}\right).
\]
Upon rationalizing, we conclude the proof of Theorem \ref{thm3}.
\end{proof}

\noindent
By making necessary substitutions, we generalize Theorem \ref{thm3} in what follows.
\begin{corollary}Let $a, b, c \in \C$, where $\Re{a}, \Re{b} \geq 0$, $\Re{c} > 0$. Then \label{ccvthm3}
\begin{equation}\label{cceqthm2}
\int_0^\infty \frac{\sin\left(bx \sqrt{x^2 + c^2}\right)}{x\sqrt{x^2 + c^2}} \eul^{-ax^2} \, \mathrm{d}x = \frac{\pi}{2c}\erf\left(\sqrt{\frac{\sqrt{a^2 c^4 + b^2 c^4} - ac^2}{2}}\right).
\end{equation}
\end{corollary}

\noindent
We do not need need to evaluate partial derivatives in \eqref{cceqthm2} due to the relationship of the equation with \eqref{eqthm1}. However, one can perform partial derivatives with respect to $a$ in \eqref{cceqthm2} and obtain an equation similar to that of \eqref{corr3}

\begin{corollary}
Let $v, z \in \mathbb{C}$, where $\Re{v} > 0$ and $\Re(v^2) > \Re(z^2)$. Then \label{corthma1}
\begin{equation}\label{eqthr1}
\erf{z} = \frac{2}{\pi}\int_0^\infty \frac{\sin\left(2vz x \sqrt{x^2 + 1}\right)}{x\sqrt{x^2 + 1}} \eul^{-(v^2 - z^2)x^2} \, \mathrm{d}x.
\end{equation}
\end{corollary}

\begin{proof}
Replacing $a$ with $v^2 - z^2$ and $b$ with $2vz$ in Theorem \ref{thm3}, we conclude the proof of Corollary \ref{corthma1}.
\end{proof}

\noindent
Theorems \ref{thm1} and \ref{thm3} lead to the following result.
\begin{corollary} Let $a, b \in \C$, where $\Re{a} > 0$, $\Re{b} \geq 0$. Then \label{corerf}
\[\int_0^b \frac{\sqrt{\sqrt{a^2 + u^2} + a}}{\sqrt{a^2 + u^2}}\eul^{\frac{a - \sqrt{a^2 + u^2}}{2}} \, \mathrm{d}u = \sqrt{2\pi} \erf\left(\sqrt{\frac{\sqrt{a^2 + b^2} - a}{2}}\right).\]
\end{corollary}

\begin{proof}
Replacing $b$ with $u$ in \eqref{eqthm1}, integrating both sides of the resulting equation with respect to $u$, from $u=0$ to $b$, we obtain
\begin{equation}\label{eqthmr1}
\int_0^\infty \frac{\sin\left(bx \sqrt{x^2 + 1}\right)}{x \sqrt{x^2 + 1}} \eul^{-ax^2} \, \mathrm{d}x = \frac{1}{2}\sqrt{\frac{\pi}{2}} \int_0^b \frac{\sqrt{\sqrt{a^2 + u^2} + a}}{\sqrt{a^2 + u^2}}\eul^{\frac{a - \sqrt{a^2 + u^2}}{2}} \, \mathrm{d}u.
\end{equation}
Equating \eqref{eqthm2} and \eqref{eqthmr1}, the proof of Corollary \ref{corerf} is complete.
\end{proof}

\noindent
From Corollary \ref{corerf}, we deduce a new integral representation for the error function, as follows.
\begin{corollary} Let $v, z \in \mathbb{C}$, where $\Re{v} > 0$ and $\Re(v^2) > \Re(z^2)$. Then \label{erfff}
\[\erf{z} = \frac{1}{\sqrt{2\pi}} \int_0^{2vz}  \frac{\sqrt{\sqrt{(v^2 - z^2)^2 + u^2} + (v^2 - z^2)}}{\sqrt{(v^2 - z^2)^2 + u^2}}\exp\left(\frac{1}{2} (v^2 - z^2) - \sqrt{(v^2 - z^2)^2 + u^2}\right) \, \mathrm{d}u.\]
\end{corollary}

\begin{proof}
Replacing $a$ with $v^2 - z^2$ and $b$ with $2vz$ in Corollary \ref{corerf}, we conclude the proof of Corollary \ref{erfff}.
\end{proof}

\begin{theorem}Let $a, b, c \in \C$, where $\Re{a}, \Re{b}, \Re{c} > 0$. Then \label{thm4}
\[\int_0^\infty \frac{\sqrt{\sqrt{x^2 + a^2} + a}}{x\sqrt{x^2 + a^2}} \sin{(cx)} \eul^{-b\sqrt{x^2 + a^2}} \, \mathrm{d}x = \frac{\pi}{\sqrt{2a}} \eul^{-ab}\erf\left(\sqrt{a\left(\sqrt{b^2 + c^2} - b\right)}\right).\]
\end{theorem}

\begin{proof}
Replacing $c$ with $u$ in \eqref{leqn2}, integrating both sides of \eqref{leqn2} with respect to $u$, from $u=0$ to $c$, we obtain
\begin{equation}\label{imp1}
\int_0^c F(a, b, u) \, \mathrm{d}u = \frac{2}{\sqrt{2a}}\eul^{-ab}\int_0^\infty \frac{\sin\left(2acw\sqrt{w^2 +1}\right)}{w\sqrt{w^2 + 1}} \eul^{-2abw^2}\, \mathrm{d}w.
\end{equation}
Utilizing the integral representation of $F(a, b, u)$ from Theorem \ref{thm2} in the left-hand side of \eqref{imp1}, and  applying Theorem \ref{thm3} to derive a closed form for the integral on the right-hand side of \eqref{imp1}, we conclude the proof of Theorem \ref{thm4}.
\end{proof}

\begin{corollary} Let $a, v, z \in \mathbb{C}$, where $\Re{a} > 0$, $\Re{v} > 0$, and $\Re\left(2av^2\right)  > \Re(z^2)$. Then  \label{corinh1}
\[\erf{z} =  \frac{\sqrt{2a}}{\pi} \eul^{\frac{2av^2 - z^2}{2}} \int_0^\infty \frac{\sqrt{\sqrt{x^2 + a^2} + a}}{x\sqrt{x^2 + a^2}} \sin\left(\sqrt{\frac{2}{a}}vz x\right) \eul^{-\frac{2av^2 - z^2}{2a}\sqrt{x^2 + a^2}} \, \mathrm{d}x.\]
\end{corollary}

\begin{proof}
Replacing $b$ with $v^2 - \frac{z^2}{2a}$ and $c$ with $\sqrt{\frac{2}{a}}vz$ in Theorem \ref{thm4}, we conclude the proof of Corollary \ref{corinh1}.
\end{proof}

\begin{corollary} Let $a, b, c \in \C$, where $\Re{a}, \Re{b}, \Re{c} > 0$. Then \label{corerf2}
\[\int_0^c \frac{\sqrt{b + \sqrt{b^2 + u^2}}}{\sqrt{b^2 + u^2}} \eul^{-a\sqrt{b^2 + u^2}} \, \mathrm{d}u = \sqrt{\frac{\pi}{a}} \eul^{-ab}\erf\left(\sqrt{a\left(\sqrt{b^2 + c^2} - b\right)}\right).\]
\end{corollary}

\begin{proof}
Utilizing the closed from of $F(a, b, u)$ from Theorem \ref{thm2} in the left-hand side of \eqref{imp1}, we deduce
\begin{equation}\label{imp2}
 \sqrt{\frac{\pi}{2}} \int_0^c \frac{\sqrt{b + \sqrt{b^2 + u^2}}}{\sqrt{b^2 + u^2}} \eul^{-a\sqrt{b^2 + u^2}} \, \mathrm{d}u = \frac{\pi}{\sqrt{2a}} \eul^{-ab}\erf\left(\sqrt{a\left(\sqrt{b^2 + c^2} - b\right)}\right).
\end{equation}
Upon simplifying \eqref{imp2}, we conclude the proof of Corollary \ref{corerf2}.
\end{proof}

\begin{corollary} Let $a, v, z \in \mathbb{C}$, where $\Re{a} > 0$, $\Re{v} > 0$, and $\Re\left(2av^2\right)  > \Re(z^2)$. Then \label{erfff2}
\[\erf{z} = \sqrt{\frac{2}{\pi}} a\eul^{\frac{2av^2 - z^2}{2}}\int_0^{\sqrt{\frac{2}{a}}vz} \frac{\sqrt{2av^2 - z^2 + \sqrt{(2av^2 - z^2)^2 + 4a^2u^2}}}{\sqrt{(2av^2 - z^2)^2 + 4a^2u^2}} \eul^{-\frac{1}{2}\sqrt{(2av^2 - z^2)^2 + 4a^2u^2}} \, \mathrm{d}u.\]
\end{corollary}

\begin{proof}
Replacing $b$ with $v^2 - \frac{z^2}{2a}$ and $c$ with $\sqrt{\frac{2}{a}}vz$ in Corollary \ref{corerf2}, we conclude the proof of Corollary \ref{erfff2}.
\end{proof}

\begin{theorem}Let $a, b\in \C$, where $\Re{a} > 0$, $\Re{b} \geq 0$. Then \label{thm5}
\begin{equation}\label{eqthm5}
\int_0^\infty \frac{\sin\left(b\sqrt{x(x+1)}\right)}{\sqrt{x+1}} \eul^{-ax} \, \mathrm{d}x = \frac{\sqrt{\sqrt{a^2 + b^2} - a}}{\sqrt{a^2 + b^2}} \eul^{\frac{a - \sqrt{a^2 + b^2}}{2}}\sqrt{\frac{\pi}{2}}.
\end{equation}
\end{theorem}

\begin{proof}
Consider the integral
\[f(a, b) =  \int_0^\infty \frac{\sin\left(ab\sqrt{x(x+1)}\right)}{\sqrt{x+1}} \eul^{-ax} \, \mathrm{d}x.\]
Evaluating the Laplace transform of $f(a, b)$ with respect to $a$, we obtain
\begin{align*}
F(s) &= \int_0^\infty \int_0^\infty \frac{\sin\left(ab\sqrt{x(x+1)}\right)}{\sqrt{x+1}} \eul^{-a(s + x)} \, \mathrm{d}a \, \mathrm{d}x = \int_0^\infty \frac{1}{\sqrt{x+1}} \frac{b\sqrt{x(x+1)}}{(s + x)^2 + b^2x(x+1)} \, \mathrm{d}x 
\\&= 2ba^{-\frac{1}{2}} (b^2 + 1)^{-\frac{3}{4}} \int_0^\infty \frac{x^2}{x^4 + \frac{b^2 + 2a}{a\sqrt{b^2 + 1}} x^2 + 1} \, \mathrm{d}x
\\&= 2ba^{-\frac{1}{2}} (b^2 + 1)^{-\frac{3}{4}} \cdot \frac{\pi}{2\sqrt{2}} \left(\frac{b^2 + 2a}{2a\sqrt{b^2 + 1}} + 1\right)^{-\frac{1}{2}}
\\&= \frac{b\pi}{\sqrt{2}} \frac{\left(1 + \sqrt{b^2 + 1}\right)^{-\frac{1}{2}}}{\sqrt{b^2 + 1}} \left(\frac{\sqrt{b^2 + 1} - 1}{2} + a\right)^{-\frac{1}{2}}.
\end{align*}
Evaluating the inverse Laplace transform of $F(s)$, we deduce
\[
f(a, b) = \sqrt{\frac{\pi}{2a}} \frac{b}{\sqrt{1 + b^2} \sqrt{1 + \sqrt{1 + b^2}}} \eul^{\frac{a - a\sqrt{1+b^2}}{2}}.
\]
Replacing $b$ with $\frac{b}{a}$ in $f(a, b)$ and simplifying further, we conclude the proof of Theorem \ref{thm5}.
\end{proof}

\begin{ex}The following result is valid:
\[
\int_0^\infty \frac{\sin\left(4\sqrt{x(x+1)}\right)}{\sqrt{x+1}} \eul^{-3x} \, \mathrm{d}x = \frac{\sqrt{\pi}}{5\eul}.\]
\end{ex}

\noindent
By differentiating both sides of \eqref{eqthm5} with respect to $b$, we derive the following corollary.
\begin{corollary} Let $a, b\in \C$, where $\Re{a} > 0$, $\Re{b} \geq 0$. Then \label{cornm}
\begin{equation}\label{eqthmm6}
\begin{split}
&\int_0^\infty \sqrt{x} \cos\left(b\sqrt{x(x+1)}\right) \eul^{-ax} \, \mathrm{d}x 
\\&\qquad= \frac{\left((a-1)\sqrt{a^2+b^2} + 2a -a^2-b^2\right)\sqrt{\sqrt{a^2+b^2}+a}}{2 \left(a^2+b^2\right)^{\frac{3}{2}}} \eul^{\frac{a-\sqrt{a^2+b^2}}{2}} \sqrt{\frac{\pi}{2}}.
\end{split}
\end{equation}
\end{corollary}

\begin{ex}The following result is valid:
\[
\int_0^\infty \sqrt{x} \cos\left(4\sqrt{x(x+1)}\right) \eul^{-3x} \, \mathrm{d}x = -\frac{9 \sqrt{\pi}}{125 \eul}.\]
\end{ex}

\noindent
By differentiating both sides of \eqref{eqthm5} twice with respect to $b$, we derive the following corollary.
\begin{corollary} Let $a, b\in \C$, where $\Re{a} > 0$, $\Re{b} \geq 0$. Then \label{cornmd2}
\begin{equation}\label{eqthmm6}
\begin{split}
&\int_0^\infty x \sqrt{x+1} \sin\left(b\sqrt{x(x+1)}\right) \eul^{-ax} \, \mathrm{d}x 
\\&\qquad= \frac{\sqrt{\sqrt{a^2+b^2}-a}}{4\left(a^2+b^2\right)^{\frac{5}{2}}} \left(9 a^2 + 2 a^3 - 3 b^2 + 2 a b^2 - a^2 b^2 - b^4 \right.
\\&\qquad\qquad\qquad\qquad \qquad\qquad\left.+ \left(6a + 4a^2 - 2b^2\right)\sqrt{a^2+b^2}\right)\eul^{\frac{a-\sqrt{a^2+b^2}}{2}} \sqrt{\frac{\pi}{2}}.
\end{split}
\end{equation}
\end{corollary}

\begin{ex}The following result is valid:
\[\int_0^\infty x \sqrt{x+1} \sin\left(4\sqrt{x(x+1)}\right) \eul^{-3x} \, \mathrm{d}x = -\frac{107 \sqrt{\pi}}{12500 \eul}.\]
\end{ex}

\noindent
By making necessary substitutions and differentiating both sides of \eqref{eqthm5} $n$ times with respect to $b$, we obtain the following result.
\begin{corollary}Let $a, b, c\in \C$, where $\Re{a}, \Re{c} > 0$, $\Re{b} \geq 0$, and $n$ a nonnegative integer. Then \label{cthm5}
\begin{align*}\label{cceqthm5}
&\int_0^\infty x^{\frac{n}{2}} (x + c)^{\frac{n - 1}{2}} \sin\left(b\sqrt{x(x+c)} + \frac{n\pi}{2}\right) \eul^{-ax} \, \mathrm{d}x 
\\&\qquad \qquad \quad = \sqrt{\frac{\pi}{2}} \frac{\partial^n}{\partial b^n}\left(\frac{\sqrt{c} \sqrt{\sqrt{a^2 c^2 + b^2 c^2} - a c}}{\sqrt{a^2 c^2 + b^2 c^2}} \eul^{\frac{a c - \sqrt{a^2 c^2 + b^2 c^2}}{2}}\right).
\end{align*}
\end{corollary}

\noindent
By differentiating both sides of \eqref{eqthm5} with respect to $a$, we derive the following corollary.
\begin{corollary} Let $a, b\in \C$, where $\Re{a}, \Re{b} > 0$. Then \label{ay1}
\begin{equation}\label{eqnha1}
\begin{split}
&\int_0^\infty \frac{x\sin\left(b\sqrt{x(x+1)}\right)}{\sqrt{x+1}} \eul^{-ax} \, \mathrm{d}x 
\\&\qquad= \frac{\left((a+1)\sqrt{a^2+b^2} +2a-a^2-b^2\right)\sqrt{\sqrt{a^2+b^2}-a}}{2 \left(a^2+b^2\right)^{\frac{3}{2}}} \eul^{\frac{a-\sqrt{a^2+b^2}}{2}} \sqrt{\frac{\pi}{2}}.
\end{split}
\end{equation}
\end{corollary}

\begin{ex}The following result is valid:
\[
\int_0^\infty \frac{x\sin\left(4\sqrt{x(x+1)}\right)}{\sqrt{x+1}} \eul^{-3x} \, \mathrm{d}x = \frac{\sqrt{\pi}}{250\eul}.\]
\end{ex}

\noindent
By differentiating both sides of \eqref{eqthm5} twice with respect to $a$, we derive the following corollary.
\begin{corollary} Let $a, b\in \C$, where $\Re{a} > 0$, $\Re{b} \geq 0$. Then \label{ay2}
\begin{equation}\label{eqthmm6}
\begin{split}
&\int_0^\infty \frac{x^2\sin\left(b\sqrt{x(x+1)}\right)}{\sqrt{x+1}} \eul^{-ax} \, \mathrm{d}x 
\\&\qquad= \frac{\sqrt{\sqrt{a^2+b^2}-a}}{4\left(a^2+b^2\right)^{\frac{5}{2}}} \left(9 a^2 - 2 a^3 + 2 a^4 - 3 b^2 - 2 a b^2 + 3 a^2 b^2 + b^4\right.
\\&\qquad\qquad\qquad\qquad \qquad\qquad\left.+ \left(6a + 2a^2 - 2a^3 - 4b^2 -2ab^2\right)\sqrt{a^2+b^2}\right)\eul^{\frac{a-\sqrt{a^2+b^2}}{2}} \sqrt{\frac{\pi}{2}}.
\end{split}
\end{equation}
\end{corollary}

\begin{ex}The following result is valid:
\[\int_0^\infty \frac{x^2\sin\left(4\sqrt{x(x+1)}\right)}{\sqrt{x+1}} \eul^{-3x} \, \mathrm{d}x = -\frac{157 \sqrt{\pi}}{12500 \eul}.\]
\end{ex}

\noindent
By making necessary substitutions and differentiating both sides of \eqref{eqthm5} $n$ times with respect to $a$, we obtain the following result.
\begin{corollary}Let $a, b, c\in \C$, where $\Re{a}, \Re{c} > 0$, $\Re{b} \geq 0$, and $n$ a nonnegative integer. Then \label{ccthm5}
\begin{align*}\label{cceqthm5}
&\int_0^\infty \frac{x^n\sin\left(b\sqrt{x(x+c)}\right)}{\sqrt{x+c}} \eul^{-ax} \, \mathrm{d}x 
\\&\qquad \qquad \quad = (-1)^n \sqrt{\frac{\pi}{2}} \frac{\partial^n}{\partial a^n}\left(\frac{\sqrt{c} \sqrt{\sqrt{a^2 c^2 + b^2 c^2} - a c}}{\sqrt{a^2 c^2 + b^2 c^2}} \eul^{\frac{a c - \sqrt{a^2 c^2 + b^2 c^2}}{2}}\right).
\end{align*}
\end{corollary}

\begin{theorem}Let $a, b, c \in \C$, where $\Re{a}, \Re{b} > 0$, $\Re{c} \geq 0$. Then \label{thm6}
\begin{equation} \label{yind2}
\int_0^\infty \frac{\sqrt{\sqrt{x^2 + a^2} - a}}{\sqrt{x^2 + a^2}} \sin(c x) \eul^{-b\sqrt{x^2 + a^2}} \, \mathrm{d}x = \frac{\sqrt{\sqrt{b^2 + c^2} - b}}{\sqrt{b^2 + c^2}} \eul^{-a\sqrt{b^2 + c^2}} \sqrt{\frac{\pi}{2}}.
\end{equation}
\end{theorem}

\begin{proof}
Let us denote
\[R(a, b, c) = \int_0^\infty \frac{\sqrt{\sqrt{x^2 + a^2} - a}}{\sqrt{x^2 + a^2}} \sin(c x) \eul^{-b\sqrt{x^2 + a^2}} \, \mathrm{d}x.\]
By substituting $x = au$, we obtain
\[R(a, b, c) = \sqrt{a} \int_0^\infty \frac{\sqrt{\sqrt{u^2 + 1} - 1}}{\sqrt{u^2 + 1}} \sin(acu) \eul^{-ab\sqrt{u^2 + 1}} \, \mathrm{d}u.\]
Further substitution of $u = \frac{\sqrt{1 -v^2}}{v}$ gives
\[R(a, b, c) = \sqrt{a} \int_0^1 v^{-\frac{3}{2}} (1 + v)^{-\frac{1}{2}} \eul^{-\frac{ab}{v}} \sin\left(\frac{a c \sqrt{1-v^2}}{v}\right) \mathrm{d}v.\]
Finally, substituting $v = \frac{1}{2w + 1}$ yields
\begin{equation}\label{leqn2m}
R(a, b, c) = \sqrt{2a}\eul^{-ab}\int_0^\infty (w + 1)^{-\frac{1}{2}} \sin\left(2ac\sqrt{w(w +1)}\right)\eul^{-2abw}\, \mathrm{d}w.
\end{equation}
Utilizing Theorem \ref{thm5} in \eqref{leqn2m}, we conclude the proof of Theorem \ref{thm6}.
\end{proof}

\begin{corollary}Let $\alpha, \beta \in \C$, where $\Re(\alpha) > 0$ and $\Re(\alpha^2) > \Re(\beta^2)$. Then \label{corfh1}
\[\int_0^\infty \frac{\sqrt{\sqrt{x^2 + \alpha^2} - \alpha}}{\sqrt{x^2 + \alpha^2}} \sin(2\alpha \beta x) \eul^{-(\alpha^2 - \beta^2)\sqrt{x^2 + \alpha^2}} \, \mathrm{d}x = \frac{\beta}{\alpha^2 + \beta^2} \eul^{-\alpha(\alpha^2 + \beta^2)} \sqrt{\pi}.\]
\end{corollary}

\begin{proof}
Substituting $b = \alpha^2 - \beta^2$, $c = 2\alpha \beta$, and $a=\alpha$ in Theorem \ref{thm6}, we conclude the proof of Corollary \ref{corfh1}.
\end{proof}

\begin{ex}The following result is valid: 
\[
\int_0^\infty \frac{\sqrt{\sqrt{x^2 + 4} - 2}}{\sqrt{x^2 + 4}} \sin(4x) \eul^{-3\sqrt{x^2 + 4}} \, \mathrm{d}x  = \frac{\sqrt{\pi}}{5 \eul^{10}}.\]
\end{ex}

\noindent
By differentiating both sides of  \eqref{yind2} with repsect to $b$, we obtain the following result.
\begin{corollary}Let $a, b, c \in \C$, where $\Re{a}, \Re{b} > 0$, $\Re{c} \geq 0$. Then \label{corfrmthm6}
\begin{align*}
\int_0^\infty &\sqrt{\sqrt{x^2 + a^2} - a}\,  \sin(c x) \eul^{-b\sqrt{x^2 + a^2}} \, \mathrm{d}x 
\\&\qquad= \frac{\sqrt{\sqrt{b^2+c^2}-b} \left((2ab + 1)\sqrt{b^2+c^2} + 2b\right)}{2 \left(b^2+c^2\right)^{\frac{3}{2}}} \eul^{-a \sqrt{b^2+c^2}} \sqrt{\frac{\pi}{2}}.
\end{align*}
\end{corollary}

\begin{ex}The following result is valid: 
\[
\int_0^\infty \sqrt{\sqrt{x^2 + 1} - 1}\, \sin(4x) \eul^{-3\sqrt{x^2 + 1}} \, \mathrm{d}x  = \frac{41\sqrt{\pi}}{250 \eul^{5}}.\]
\end{ex}

\noindent
By differentiating both sides of  \eqref{yind2} twice with repsect to $b$, we obtain the following result.
\begin{corollary}Let $a, b, c \in \C$, where $\Re{a}, \Re{b} > 0$, $\Re{c} \geq 0$. Then \label{cor2fromthm6}
\begin{align*}
\int_0^\infty &\sqrt{\sqrt{x^2 + a^2} - a} \sqrt{x^2 + a^2}  \sin(c x) \eul^{-b\sqrt{x^2 + a^2}} \, \mathrm{d}x 
\\&\qquad= \frac{\sqrt{\sqrt{b^2+c^2}-b}}{4 \left(b^2+c^2\right)^{\frac{5}{2}}} \left(9 b^2 + 4 a b^3 + 4 a^2 b^4 - 3 c^2 + 4 a b c^2 + 4 a^2 b^2 c^2\right.
\\&\qquad \qquad \qquad \qquad \qquad \quad \left.  + (8 a b^2 + 6 b - 4 a c^2 ) \sqrt{b^2 + c^2}\right) \eul^{-a \sqrt{b^2+c^2}} \sqrt{\frac{\pi }{2}}.
\end{align*}
\end{corollary}

\begin{ex}The following result is valid: 
\[
\int_0^\infty \sqrt{\sqrt{x^2 + 1} - 1}\sqrt{x^2 + 1}\, \sin(4x) \eul^{-3\sqrt{x^2 + 1}} \, \mathrm{d}x  = \frac{1363 \sqrt{\pi}}{12500\eul^{5}}.\]
\end{ex}

\noindent
By differentiating both sides of  \eqref{yind2} with repsect to $c$, we obtain the following result.
\begin{corollary}Let $a, b, c \in \C$, $\Re{a}, \Re{b} > 0$, $\Re{c} \geq 0$. Then \label{yinde3}
\begin{align*}
\int_0^\infty &\frac{x\sqrt{\sqrt{x^2 + a^2} - a}}{\sqrt{x^2 + a^2}}  \cos(c x) \eul^{-b\sqrt{x^2 + a^2}} \, \mathrm{d}x 
\\&\qquad= \frac{\sqrt{\sqrt{b^2+c^2}-b} \left(\left(b - 2ac^2\right)\sqrt{b^2 + c^2} + b^2 - c^2 \right)}{2c\left(b^2+c^2\right)^{\frac{3}{2}}} \eul^{-a \sqrt{b^2+c^2}} \sqrt{\frac{\pi}{2}}.
\end{align*}
\end{corollary}

\begin{ex}The following result is valid: 
\[
\int_0^\infty \frac{x\sqrt{\sqrt{x^2 + 1} - 1}}{\sqrt{x^2 + 1}} \cos(4x) \eul^{-3\sqrt{x^2 + 1}} \, \mathrm{d}x  = -\frac{19\sqrt{\pi}}{125\eul^{5}}.\]
\end{ex}

\noindent
By differentiating both sides of  \eqref{yind2} twice with repsect to $c$, we obtain the following result.
\begin{corollary}Let $a, b, c \in \C$, $\Re{a}, \Re{b} > 0$, $\Re{c} \geq 0$. Then \label{yinde4}
\begin{align*}
\int_0^\infty &\frac{x^2\sqrt{\sqrt{x^2 + a^2} - a}}{\sqrt{x^2 + a^2}}  \sin(c x) \eul^{-b\sqrt{x^2 + a^2}} \, \mathrm{d}x 
\\&\qquad= \frac{\sqrt{\sqrt{b^2+c^2}-b}}{4 \left(b^2+c^2\right)^{\frac{5}{2}}} \left(9 b^2 + 4 a b^3 - 3 c^2 + 4 a b c^2 - 4 a^2 b^2 c^2 - 4 a^2 c^4\right.
\\&\qquad \qquad \qquad \qquad \qquad \quad \left.  + (8 a b^2 + 6 b - 4 a c^2 ) \sqrt{b^2 + c^2}\right) \eul^{-a \sqrt{b^2+c^2}} \sqrt{\frac{\pi }{2}}.
\end{align*}
\end{corollary}

\begin{ex}The following result is valid: 
\[
\int_0^\infty \frac{x^2\sqrt{\sqrt{x^2 + 1} - 1}}{\sqrt{x^2 + 1}} \sin(4x) \eul^{-3\sqrt{x^2 + 1}} \, \mathrm{d}x  = -\frac{1137\sqrt{\pi}}{12500\eul^{5}}.\]
\end{ex}

\noindent
By differentiating both sides of  \eqref{yind2} with respect to $a$, we obtain the following result.
\begin{corollary}Let $a, b, c \in \C$, $\Re{a}, \Re{b} > 0$, $\Re{c} \geq 0$. Then \label{yinde5}
\begin{align*}
\int_0^\infty &\frac{\sqrt{\sqrt{x^2+a^2}-a} \left((2ab + 1) \sqrt{x^2+a^2} + 2a\right) \sin(cx)}{\left(x^2+a^2\right)^{\frac{3}{2}}} \eul^{-b \sqrt{x^2+a^2}} \, \mathrm{d}x 
\\&\qquad\quad= \sqrt{\sqrt{b^2 + c^2} - b} \eul^{-a \sqrt{b^2+c^2}} \sqrt{2\pi}.
\end{align*}
\end{corollary}

\begin{ex}The following result is valid: 
\[
\int_0^\infty \frac{\sqrt{\sqrt{x^2+1}-1} \left(7\sqrt{x^2+ 1} + 2\right) \sin(4x)}{\left(x^2+1\right)^{\frac{3}{2}}} \eul^{-3\sqrt{x^2+1}} \, \mathrm{d}x  = \frac{2\sqrt{\pi}}{\eul^{5}}.\]
\end{ex}

\noindent
By differentiating both sides of  \eqref{yind2} twice with repsect to $a$, we obtain the following result.
\begin{corollary}Let $a, b, c \in \C$, $\Re{a}, \Re{b} > 0$, $\Re{c} \geq 0$. Then \label{yinde6}
\begin{align*}
\int_0^\infty &\frac{\sqrt{\sqrt{x^2+a^2}-a}}{\left(x^2+a^2\right)^{\frac{5}{2}}} \left(9 a^2 + 4 a^3 b + 4 a^4 b^2 - 3 x^2 + 4 a b x^2 + 4 a^2 b^2 x^2 \right.
\\& \qquad \qquad \qquad \qquad \left.+ \left(8a^2b + 6a - 4bx^2\right) \sqrt{x^2+a^2}\right) \sin(cx) \eul^{-b \sqrt{x^2+a^2}} \, \mathrm{d}x 
\\&\qquad= 2\sqrt{\sqrt{b^2 + c^2} - b} \sqrt{b^2+c^2} \eul^{-a \sqrt{b^2+c^2}} \sqrt{2\pi}.
\end{align*}
\end{corollary}

\begin{ex}The following result is valid: 
\[\int_0^\infty \frac{\sqrt{\sqrt{x^2+1}-1} \left(\left(30-12 x^2\right) \sqrt{x^2+1} + 57 + 45 x^2\right) \sin(4x)}{\left(x^2+1\right)^{\frac{5}{2}}} \eul^{-3\sqrt{x^2+1}} \, \mathrm{d}x  = \frac{20\sqrt{\pi}}{\eul^{5}}.\]
\end{ex}

\subsection{More Identities} \label{sec2.2}
In this subsection, we use a formula from Cauchy's residue theorem to derive an identity for the inverse Laplace transform. By applying this identity, along with a result from Srivastava and Y\"{u}rekli's work, we present additional new and interesting identities.
\begin{lemma}[Cauchy]\label{thm2.7h}
Let $g(x)$ be a continuous and sufficiently smooth function defined on $[0, \infty)$. Suppose that $g(z)$ is the extension of $g(x)$ into the complex plane and is analytic in the upper half-plane. Then for $\Re(\lambda) > 0$, the following holds:
\[
\int_0^\infty \frac{\Re{g(\img x)}}{x^2 + \lambda^2} \, \mathrm{d}x = \frac{\pi g(\lambda)}{2\lambda}, 
\]
provided the integral exists.
\end{lemma}

\begin{proof}
We first rewrite $\Re(g(ix))$ using the identity 
\[\Re(g(ix)) = \frac{g(ix) + g(-ix)}{2} = \Re(g(-ix)).\]
Thus, the integral becomes
\[\int_0^{\infty} \frac{\Re(g(ix))}{ x^2 + \lambda^2} \, \mathrm{d}x = \frac{1}{2}\int_{-\infty}^{\infty}  \frac{\Re(g(-ix))}{x^2 + \lambda^2} \, \mathrm{d}x.\]
Considering the contour integral along the real axis and a semi-circular arc in the upper half-plane, and applying Cauchy's residue theorem \cite[\S1.10(iv)]{bib23}, the result follows.
\end{proof}
\noindent
In the next theorem, we provide a formula in terms of the inverse Laplace transform.
\begin{lemma} \label{thmh8}
Let $g$ be as stated in Lemma \ref{thm2.7h}. Then for $ \Re{b} > 0$,
\[\int_0^\infty \frac{\Re{g(\img x)}}{x} \sin(bx) \, \mathrm{d}x = \frac{\pi}{2}\mathcal{L}^{-1}\left\{\frac{g(s)}{s}\right\}(b),\]
provided the integral exists.
\end{lemma}

\begin{proof}
From Lemma \ref{thm2.7h}, we derive for $\Re{s} > 0$,
\[\int_0^\infty \eul^{-st} \left(\int_0^\infty \frac{\Re{g(\img x)}}{x} \sin(btx) \, \mathrm{d}t\right) \mathrm{d}x = \frac{\pi g\left(\frac{s}{b}\right)}{2s}.\]
Interchanging the order of integration, which is permissible by the absolute convergence of the integrals involved, we have
\[\int_0^\infty \eul^{-st} \left(\int_0^\infty \frac{\Re{g(\img x)}}{x} \sin(btx) \, \mathrm{d}x\right) \mathrm{d}t = \frac{\pi g\left(\frac{s}{b}\right)}{2s}.\]
Evaluating the inverse Laplace transform, we deduce
\[\int_0^\infty \frac{\Re{g(\img x)}}{x} \sin(bx) \, \mathrm{d}x =\frac{\pi}{2} \mathcal{L}^{-1}\left\{\frac{g\left(\frac{s}{b}\right)}{s}\right\}(t) \biggr\vert_{t=1}= \frac{\pi}{2}\mathcal{L}^{-1}\left\{\frac{g\left(\frac{s}{b}\right)}{s}\right\}(1).\]
It now remains to show that
\[\mathcal{L}^{-1}\left\{\frac{g\left(\frac{s}{b}\right)}{s}\right\}(1) = \mathcal{L}^{-1}\left\{\frac{g(s)}{s}\right\}(b).\]
Denote $\phi(t) = \mathcal{L}^{-1}\left\{\frac{g\left(\frac{s}{b}\right)}{s}\right\}(t)$. This implies
\[\frac{g\left(\frac{s}{b}\right)}{s} = \int_0^\infty \phi(t) \eul^{-st} \, \mathrm{d}t.\]
Replacing $s$ with $sb$, and making the substitution $t \rightarrow tb$, we derive
\[\frac{g(s)}{sb} = \frac{1}{b} \int_0^\infty \phi\left(\frac{t}{b}\right) \eul^{-st} \, \mathrm{d}t.\]
Therefore,
\[\phi\left(\frac{t}{b}\right) = \mathcal{L}^{-1}\left\{\frac{g(s)}{s}\right\}(t).\]
On the substitution $t \rightarrow tb$, we derive
\[\phi(t) = \mathcal{L}^{-1}\left\{\frac{g(s)}{s}\right\}(tb).\]
Setting $t=1$, we deduce
\[\phi(1) = \mathcal{L}^{-1}\left\{\frac{g(s)}{s}\right\}(b) = \mathcal{L}^{-1}\left\{\frac{g\left(\frac{s}{b}\right)}{s}\right\}(1).\]
This completes the proof of Lemma \ref{thmh8}.
\end{proof}

\begin{corollary} Let $a, b \in \C$, where $\Re{a}, \, \Re{b} > 0$. Then \label{newcoy}
\[\int_0^\infty \frac{\sqrt{\sqrt{x^2 + a^2} + a}}{x} \sin{(bx)} \, \mathrm{d}x = \sqrt{\frac{\pi}{2b}} \eul^{-ab} + \sqrt{\frac{a}{2}} \pi \erf\left(\sqrt{ab}\right).\]
\end{corollary}

\begin{proof}
Utilizing $g(x) = \sqrt{x + a}$ in Lemma \ref{thmh8}, we conclude the proof of Corollary \ref{newcoy}.
\end{proof}

\begin{corollary} Let $a, b \in \C$, where $\Re{a} > 0, \, \Re{b} \geq 0$. Then \label{newcjk1}
\[\int_0^\infty \frac{\sqrt{\sqrt{x^2 + a^2} + a}}{x\sqrt{x^2 + a^2}} \sin{(bx)} \, \mathrm{d}x =  \frac{\pi}{\sqrt{2a}} \erf\left(\sqrt{ab}\right).\]
\end{corollary}

\begin{proof}
Utilizing $g(x) = \frac{1}{\sqrt{x + a}}$ in Lemma \ref{thmh8}, we conclude the proof of Corollary \ref{newcjk1}.
\end{proof}

\begin{theorem}\label{corfst}
Suppose that $\Im{g(\img x)}$ is absolutely and square integrable on $[0, \infty)$. Then for $\Re{\mu} > 0$,
\[\int_0^\infty \frac{\Im{g(\img x)}}{x^{\mu}} \, \mathrm{d}x = \frac{\cosec\left(\frac{\pi \mu}{2}\right)}{\Gamma(\mu)} \sqrt{\frac{\pi}{2}} \int_0^{\infty} x^{\mu -1} \mathcal{F}_s\{\Im{g(\img t)}\}(x) \, \mathrm{d}x,\]
provided that both integrals exists.
\end{theorem}

\begin{proof}
Srivastava and Y\"{u}rekli \cite[(10)]{stava} proved that
\begin{equation} \label{applhur1}
\int_0^\infty \frac{g(x)}{x^{\mu}} \, \mathrm{d}x = \frac{1}{\Gamma(\mu)}  \int_0^{\infty} x^{\mu -1} G_L(x) \, \mathrm{d}x,
\end{equation}
where $G_L(x)$ is the Laplace transform of $g(t)$. Suppose $h(x)$ also satifies the conditions on $\Im{g(\img x)}$, then \cite[(1.14.14)]{bib23}
\begin{equation}\label{feqh2}
\int_0^\infty h(x) \Im{g(\img x)} \, \mathrm{d}x = \int_0^\infty \mathcal{F}_s\{h(x)\}(a)  \, \mathcal{F}_s\{\Im{g(\img x)}\}(a)  \, \mathrm{d}a.
\end{equation}
By making the substitutions $x \to t$, $h(x) = \eul^{-xt}$ in \eqref{feqh2}, we derive 
\begin{equation} \label{applhur2}
\mathcal{L}\{\Im{g(\img t)}\}(x) = \sqrt{\frac{2}{\pi}} \int_0^\infty \frac{a}{x^2 + a^2} \mathcal{F}_s\{\Im{g(\img x)}\}(a)\, \mathrm{d}a.
\end{equation}
By making the substitution $g(t) \rightarrow \Im{g(\img t)}$ and utilizing \eqref{applhur2} in \eqref{applhur1}, we obtain
\[
\int_0^\infty \frac{\Im{g(\img x)}}{x^{\mu}} \, \mathrm{d}x = \frac{1}{\Gamma(\mu)} \int_0^{\infty} \int_0^{\infty}  \frac{ax^{\mu -1}}{x^2 + a^2}  \mathcal{F}_s\{\Im{g(\img x)}\}(a) \, \mathrm{d}a  \, \mathrm{d}x.\]
Interchanging the order of integration, which is permissible by the absolute convergence of the integrals involved, and evaluating the inner integral, we conclude the proof of Theorem \ref{corfst}.
\end{proof}

\noindent
We apply Theorem \ref{corfst} in the following corollary.
\begin{corollary} Let $a, \mu \in \C$, where $\Re{a} > 0$, $-\frac{1}{2} < \Re{\mu} < 0$. Then \label{hafeeza1}
\[\int_0^\infty \frac{\sqrt{\sqrt{x^2+a^2} + a}}{\sqrt{x^2 + a^2}} x^{-\mu-1} \, \mathrm{d}x = -\frac{2^{\frac{1}{2}-2\mu} a^{-\frac{1}{2}-\mu} \Gamma(2\mu)}{\mu\Gamma(\mu)^2} \cosec\left(\frac{\pi \mu}{2}\right) \pi.\]
\end{corollary}

\begin{proof}
Utilizing $g(x) = \frac{1}{x\sqrt{x + a}}$ in Theorem \ref{corfst}, we obtain
\begin{equation}\label{eqnhafj}
\int_0^\infty \frac{\sqrt{\sqrt{x^2+a^2} + a}}{\sqrt{x^2 + a^2}} x^{-\mu-1} \, \mathrm{d}x = \frac{\cosec\left(\frac{\pi \mu}{2}\right)}{\Gamma(\mu)} \sqrt{\frac{\pi}{2}} \int_0^\infty  x^{\mu - 1} \mathcal{F}_s\left\{\frac{\sqrt{\sqrt{t^2+a^2} + a}}{t\sqrt{t^2 + a^2}}\right\}(x) \, \mathrm{d}x.
\end{equation}
Applying Corollary \ref{newcjk1} in \eqref{eqnhafj}, we deduce
\begin{equation}\label{eqnhafj2}
\int_0^\infty \frac{\sqrt{\sqrt{x^2+a^2} + a}}{\sqrt{x^2 + a^2}} x^{-\mu-1} \, \mathrm{d}x = \frac{\cosec\left(\frac{\pi \mu}{2}\right)}{\Gamma(\mu)}  \frac{\pi}{\sqrt{2a}} \mathcal{M}\left\{\erf\left(\sqrt{ax}\right)\right\}(\mu).
\end{equation}
By evaluating the Mellin transform in \eqref{eqnhafj2}, we conclude the proof of Corollary \ref{hafeeza1}.
\end{proof}

\section{Conclusion}
Most of the formulas we have presented in this work are non-exhaustive. For instance, many of the theorems presented in Section \ref{sec2.1} can generate more new generalized integrals. These generalized integrals can either be represented in terms of the Laplace or Fourier transform. Corollaries \ref{cor2}, \ref{cor3}, \ref{cthm5}, and \ref{ccthm5} can yield infinitely many results, though complications arise when $n$ is large, as it becomes difficult to compute the derivatives on the right-hand side of these corollaries. We have only provided the first two derivatives with respect to the variables $a$ and $b$ in the aforementioned corollaries. Lemma \ref{thm2.7h} and Theorem \ref{corfst} can be employed to derive additional identities. A notable example that follows from Corollary \ref{hafeeza1} is:
\[\int_0^{\infty} \frac{\sqrt{\sqrt{x^2 + 1} + 1}}{x^{\frac34}\sqrt{x^2 + 1}} \, \mathrm{d}x = \frac{2\pi^{\frac{3}{2}}}{\sqrt{2 - \sqrt{2}} \Gamma^2\left(\frac{3}{4}\right)}.
\]

\section*{Acknowledgment}
The first author would like to express his heartfelt gratitude to the Spirit of Ramanujan (SOR) STEM Talent Initiative, directed by the Marvin Rosenblum Professor of Mathematics, Ken Ono, for providing the essential computational tools that facilitated the verification and validation of results presented in this paper. Additionally, the first author would like to thank Dr.~H.~P.~Adeyemo, Dr.~D.~A.~Dikko, Mr.~G.~S.~Lawal, and Dr.~I.~Adinya from the Department of Mathematics at the University of Ibadan for their unwavering support and encouragement.

\section*{Conflict of interest}
The authors declare that there are no relevant financial or non-financial competing interests to report in relation to this article.

\section*{Funding}
The author did not receive funding from any organization for the submitted work.

\section*{Appendix}
Here, we provide the detailed steps for evaluating $F(s)$, as referenced in the proofs of Theorems \ref{thm1} and \ref{thm3}.
\begin{enumerate}[label = \textbf{\Alph*}:]
\item \label{appendix1} By interchanging the order of integration, which is justified by the absolute convergence of the integrals involved, we derive
\begin{align*}
F(s) &= \int_0^\infty \int_0^\infty \cos\left(abx \sqrt{x^2 + 1}\right) \eul^{-a(s + x^2)} \, \mathrm{d}a\, \mathrm{d}x = \int_0^\infty \frac{s + x^2}{(s + x^2)^2 + b^2 x^2(x^2 + 1)} \, \mathrm{d}x
\\&= \int_0^\infty \frac{s + x^2}{\left(1+b^2\right)x^4 + (2s + b^2)x^2 + s^2} \, \mathrm{d}x = \frac{1}{1 + b^2} \int_0^\infty \frac{s + x^2}{x^4 + \frac{2s + b^2}{1+b^2}x^2 + \frac{s^2}{1+b^2}} \, \mathrm{d}x.
\end{align*}
Applying the quadratic formula, we find
\[
F(s) = \frac{1}{1 + b^2} \int_0^\infty \frac{s + x^2}{\left(x^2 + \frac{b^2+2s + b\sqrt{b^2-4 s^2+4 s}}{2 \left(1+b^2\right)}\right)\left(x^2 + \frac{b^2+2s - b\sqrt{b^2-4 s^2+4 s}}{2 \left(1+b^2\right)}\right)} \, \mathrm{d}x.
\]
Decomposing into partial fractions, we obtain
\begin{align*}
F(s) &= \frac{s}{b\sqrt{b^2 - 4s^2 + 4s}} \int_0^\infty \left(\frac{1}{x^2 + \frac{b^2+2s - b\sqrt{b^2-4 s^2+4 s}}{2 \left(1+b^2\right)}} - \frac{1}{x^2 + \frac{b^2+2s + b\sqrt{b^2-4 s^2+4 s}}{2 \left(1+b^2\right)}}\right) \mathrm{d}x
\\&\quad + \frac{1}{2\left(1+b^2\right)} \int_0^\infty \left(\frac{1}{x^2 + \frac{b^2+2s - b\sqrt{b^2-4 s^2+4 s}}{2 \left(1+b^2\right)}} + \frac{1}{x^2 + \frac{b^2+2s + b\sqrt{b^2-4 s^2+4 s}}{2 \left(1+b^2\right)}}\right) \mathrm{d}x
\\&\quad - \frac{b^2 + 2s}{2(1+b^2)^2}\int_0^\infty \frac{1}{\left(x^2 + \frac{b^2+2s + b\sqrt{b^2-4 s^2+4 s}}{2 \left(1+b^2\right)}\right)\left(x^2 + \frac{b^2+2s - b\sqrt{b^2-4 s^2+4 s}}{2 \left(1+b^2\right)}\right)} \, \mathrm{d}x
\\&= \left(\frac{s}{b\sqrt{b^2 - 4s^2 + 4s}}  - \frac{b^2 + 2s}{2b(1+b^2)\sqrt{b^2-4 s^2+4 s}}\right) 
\\&\qquad \times \int_0^\infty \left(\frac{1}{x^2 + \frac{b^2+2s - b\sqrt{b^2-4 s^2+4 s}}{2 \left(1+b^2\right)}} - \frac{1}{x^2 + \frac{b^2+2s + b\sqrt{b^2-4 s^2+4 s}}{2 \left(1+b^2\right)}}\right) \mathrm{d}x
\\&\quad + \frac{1}{2\left(1+b^2\right)} \int_0^\infty \left(\frac{1}{x^2 + \frac{b^2+2s - b\sqrt{b^2-4 s^2+4 s}}{2 \left(1+b^2\right)}} + \frac{1}{x^2 + \frac{b^2+2s + b\sqrt{b^2-4 s^2+4 s}}{2 \left(1+b^2\right)}}\right) \mathrm{d}x
\\&= \frac{1}{1 + b^2}\left(\frac{1}{2} -\frac{b (2 s-1)}{2 \sqrt{b^2-4s^2+4 s}}\right) \int_0^\infty \frac{\mathrm{d}x}{x^2 + \frac{b^2+2s + b\sqrt{b^2-4 s^2+4 s}}{2 \left(1+b^2\right)}}
\\&\quad+\frac{1}{1 + b^2}\left(\frac{1}{2} + \frac{b (2 s-1)}{2 \sqrt{b^2-4s^2+4 s}}\right) \int_0^\infty \frac{\mathrm{d}x}{x^2 + \frac{b^2+2s - b\sqrt{b^2-4 s^2+4 s}}{2 \left(1+b^2\right)}}.
\end{align*}
After evaluating the integrals, we derive
\begin{align*}
F(s) &= \frac{1}{1 + b^2}\left(\frac{1}{2} -\frac{b (2 s-1)}{2 \sqrt{b^2-4s^2+4 s}}\right) \frac{\pi \sqrt{2(1+b^2)}}{2 \sqrt{b^2+2s + b\sqrt{b^2-4 s^2+4 s}}}
\\&\quad+\frac{1}{1 + b^2}\left(\frac{1}{2} + \frac{b (2 s-1)}{2 \sqrt{b^2-4s^2+4 s}}\right)   \frac{\pi \sqrt{2(1+b^2)}}{2 \sqrt{b^2+2s - b\sqrt{b^2-4 s^2+4 s}}}.
\end{align*}
By rationalizing, we get
\begin{align*}
F(s) &= \frac{1}{1 + b^2}\left(\frac{1}{2} -\frac{b (2 s-1)}{2 \sqrt{b^2-4s^2+4 s}}\right) \frac{\pi \sqrt{b^2+2s - b\sqrt{b^2-4 s^2+4 s}}}{2\sqrt{2}s}
\\&\quad+\frac{1}{1 + b^2}\left(\frac{1}{2} + \frac{b (2 s-1)}{2 \sqrt{b^2-4s^2+4 s}}\right)  \frac{\pi \sqrt{b^2+2s + b\sqrt{b^2-4 s^2+4 s}}}{2\sqrt{2}s}
\\&= \frac{\pi}{4\sqrt{2}s\left(1+b^2\right)}\left(\sqrt{b^2 + 2s - b\sqrt{b^2 - 4s^2 + 4s}} + \sqrt{b^2 + 2s + b\sqrt{b^2 - 4s^2 + 4s}}  \right.
\\&\quad \left.+ \frac{b(2s-1)\left(\sqrt{b^2 + 2s + b\sqrt{b^2 - 4s^2 + 4s}} -\sqrt{b^2 + 2s - b\sqrt{b^2 - 4s^2 + 4s}}\right)}{\sqrt{b^2 - 4s^2 + 4s}}\right).
\end{align*}
Utilizng the identity $\sqrt{x} \pm \sqrt{y} = \sqrt{x + y \pm 2\sqrt{xy}}$ for 
\[\sqrt{b^2 + 2s - b\sqrt{b^2 - 4s^2 + 4s}} \pm \sqrt{b^2 + 2s + b\sqrt{b^2 - 4s^2 + 4s}},\]
we find
\[
F(s) = \frac{\pi}{4s\left(1+b^2\right)}\left(\sqrt{b^2 + 2s + 2s \sqrt{1 + b^2}} + \frac{b(2s - 1)\sqrt{b^2 + 2s - 2s \sqrt{1 + b^2}}}{\sqrt{b^2 - 4s^2 + 4s}}\right).
\]
Rationalizing the second term and simplifying, we obtain
\begin{align*}
F(s) &= \frac{\pi}{4s\left(1+b^2\right)}\left(\sqrt{b^2 + 2s + 2s \sqrt{1 + b^2}} + \frac{b^2(2s - 1)}{\sqrt{b^2 + 2s + 2s \sqrt{1 + b^2}}}\right) 
\\&= \frac{\pi}{2\left(1+b^2\right)} \frac{1 + b^2 + \sqrt{1 + b^2}}{\sqrt{b^2 + 2s + 2s \sqrt{1 + b^2}}}
\\&= \frac{\pi\left(1 + b^2 + \sqrt{1 + b^2}\right)}{2\sqrt{2}(1+b^2)\sqrt{1 + \sqrt{1 + b^2}}} \frac{1}{\left(s + \frac{b^2}{2\left(1 + \sqrt{1 + b^2}\right)}\right)^{\frac12}}.
\end{align*}

\item \label{appendix2} By interchanging the order of integration, which is permissible due to the absolute convergence of the integrals, we obtain
\begin{align*}
F(s) &= \int_0^\infty \int_0^\infty \frac{\sin\left(abx \sqrt{x^2 + 1}\right)}{x\sqrt{x^2 + 1}} \eul^{-a(s + x^2)} \, \mathrm{d}a\, \mathrm{d}x = \int_0^\infty \frac{b}{(s + x^2)^2 + b^2 x^2(x^2 + 1)} \, \mathrm{d}x
\\&= \int_0^\infty \frac{b}{\left(1+b^2\right)x^4 + (2s + b^2)x^2 + s^2} \, \mathrm{d}x = \frac{b}{1 + b^2} \int_0^\infty \frac{1}{x^4 + \frac{2s + b^2}{1+b^2}x^2 + \frac{s^2}{1+b^2}} \, \mathrm{d}x.
\end{align*}
Next, we apply the quadratic formula, yielding
\[
F(s) = \frac{b}{1 + b^2} \int_0^\infty \frac{1}{\left(x^2 + \frac{b^2+2s + b\sqrt{b^2-4 s^2+4 s}}{2 \left(1+b^2\right)}\right)\left(x^2 + \frac{b^2+2s - b\sqrt{b^2-4 s^2+4 s}}{2 \left(1+b^2\right)}\right)} \, \mathrm{d}x.
\]
Resolving into partial fractions, we have
\[
F(s) = \frac{1}{\sqrt{b^2 - 4s^2 + 4s}} \int_0^\infty \left(\frac{1}{x^2 + \frac{b^2+2s - b\sqrt{b^2-4 s^2+4 s}}{2 \left(1+b^2\right)}} - \frac{1}{x^2 + \frac{b^2+2s + b\sqrt{b^2-4 s^2+4 s}}{2 \left(1+b^2\right)}}\right) \mathrm{d}x.
\]
Evaluating the integral, we derive
\[
F(s) =  \frac{1}{\sqrt{b^2 - 4s^2 + 4s}} \left(\frac{\pi \sqrt{2(1+b^2)}}{2 \sqrt{b^2+2s + b\sqrt{b^2-4 s^2+4 s}}} -  \frac{\pi \sqrt{2(1+b^2)}}{2 \sqrt{b^2+2s - b\sqrt{b^2-4 s^2+4 s}}}\right).
\]
\end{enumerate}

\bibliographystyle{unsrt}
\bibliography{Submission_Laplace_Transform}

\end{document}